# Arithmetization of well-formed parenthesis strings. Motzkin Numbers of the Second Kind

Gennady Eremin

ergenns@gmail.com

December 23, 2020

**Abstract**. In this paper, we perform an arithmetization of well-formed parenthesis strings with zeros (Motzkin words) and of corresponding Motzkin paths. The transformations used are reminiscent of Gödel numbering for mathematical objects of some formal language. We construct a Motzkin series that is as close as possible to natural numbers by many formal features. Parenthesis strings are encoded by ternary codes and corresponding natural numbers, *Motzkin numbers of the 2nd kind*, which made it possible to formalize and simplify the analysis of the successor function, to specify and clarify the procedure for selecting a successor. In the process of arithmetization of well-wormed parenthesis strings, various special numbers appeared. In this regard, in conclusion we will talk a little about numbers of the form $3^n+2$ and "mirror" numbers $2\times3^n+1$.

*Keywords*: strong order, arithmetization of formal system, Motzkin words, Motzkin numbers, successor function, special numbers of the form $3^n+2$ and $2\times3^n+1$.

# 1 Introduction

Recently, much interest has been shown in the ordering of well-formed parenthesis strings with zeros. For example, in [Fan19, GZ14] the authors consider the partial order on Motzkin words, and in the article [BP14] a specific distance between Motzkin words in the Tamari lattice is established. In this paper, we consider a sequence of digitized Motzkin words, which is obtained in accordance with the formal features of a series of natural numbers. This article continues [Ere20].

The arithmetization method makes it possible to establish a simple one-to-one mapping of the set of words (in the language alphabet) into a subset of images, natural numbers. The image of a word is often referred to as a number. Relations and operations defined in words are transformed with such a mapping into relations and operations on natural numbers.

In this work, using the arithmetization method, we build a two-level model of balanced parentheses. The lower baseline is the foundation of such a model. At the lower level, the images of Motzkin words guarantee a bijective representation in the form of ternary codes with the preservation of external parameters and rules for constructing parenthesis strings (word length, balance of parentheses, etc.). Ternary codes simplify the procedure for strict ordering of sequence elements and allow calculating the numerical increment of the successor function (successor oper-ation, "stroke" operation). To implement the successor function, the article considers a combinatorial search algorithm (directed enumeration of options).

The second level of the parenthesis string model is a lower-level add-on. At the second level, we work with natural numbers, *Motzkin numbers of the 2nd kind*, which correspond to ternary words. In conclusion, we will talk a little about the accompanying special numbers of the form $3^n+2$ and $2\times3^n+1$, which resemble Fermat's numbers and Mersenne numbers and are little known in the literature.

Some definitions begin the natural numbers with zero since the work of a group of French mathematicians Nicolas Bourbaki defined natural numbers as the cardinalities of finite sets. It was Bourbaki who first introduced the symbol $\varnothing$ for an empty set with cardinality 0. In the natural series, 0 is important to us, since we also start Motzkin words with 0. In the alphabet of natural numbers and in the alphabet of Motzkin words, the sign 0 plays an essential role, respectively, the number 0 and the word 0.

The Motzkin word alphabet has three characters: a zero, a left parenthesis, and a right parenthesis. We will encode the left parenthesis with one and the right parenthesis with two. Obviously, no two different Motzkin words will give us the same ternary code (and corresponding decimal code). The encoding in the alphabet $A = \{0, 1, 2\}$ allows us to work both with words (you can use the basic operation of *concatenation* of words, as well as the *empty word $\varepsilon$*) and with digital codes using arithmetic operations. We denote by $A^*$ the set of finite words in the alphabet $A$ including $\varepsilon$.

In a Motzkin word, the parentheses (they are 1's and 2's) are balanced, i.e. the number of 1's and 2's is the same, and in any word *prefix* (initial fragment) the number of 1's is not less than 2's. We will work more often with word suffixes (code ends), so we will reformulate the second rule for suffixes: *in any suffix of the Motzkin code, the number of 1's is not more than 2's, otherwise the balance is disturbed.* A Motzkin word can be without zeros, or it can contain only zeros.

Motzkin words of length $n$, $n$-words, are counted by Motzkin numbers $M_n$. Here is the beginning of the sequence A001006, in which the $n$th element is equal to the number of $n$-words, $n = 0, 1, \ldots$: 1, 1, 2, 4, 9, 21, 51, 127, 323, ...

The existence of a single empty word $\varepsilon$ is allowed, so $M_0 = 1$. If possible, we will do without the empty word (note that there is no empty number among natural numbers either). Another example, in article [Fan19] the author also works without the Motzkin number $M_0$.

We construct an ordered sequence of ternary Motzkin words, focusing on the formal traits of the set of natural numbers, namely: (a) all words are unique, (b) words are sorted by code length (range distribution, *external sorting*), and finally (c) within each range, words are sorted according to the alphabetical order (*inter-* *nal sorting*). The uniqueness of Motzkin words is achieved by the following statement.

**Proposition 1.1.** *We consider two Motzkin words to be identical, if they differ only by a zero prefix.*

We do not write zeros in front of a natural number (zero prefix), since this is a ballast. The null prefix does not change an integer, for example, 12 = 012 = 0012. The situation is similar with Motzkin words, only the word 0 has a zero prefix, and according to Proposition 1.1  0 = 00 = 000 … and so on. Other ternary words begin with the symbol 1 and together with the word 0 form the set of *unique* ternary Motzkin words $M = \{w_i\} \subset A^*$.

When writing strings, we will use the concatenation operation, this is when words are connected to each other using the sign · (analogous to multiplication of variables). Like multiplication, the concatenation operation generates a degree; this is when a character or group of characters is repeated in words. For example, $1 \cdot 0^3 \cdot 2 = 10002$ or $11 \cdot (122)^2 \cdot 0 = 111221220$ (multiplication signs are often omitted). The zero degree of any word is the concatenation of zero words, and this is considered an empty word.

In ranges, we sort words in ascending order according to the alphabetical order: $0 < 1 < 2$. For example, let the words $x = pas$ and $y = pbt$ of the same length have the same prefix $p \in A^*$ (if $p$ has only one character, then $p = 1$) and let them, in general, have different suffixes $s, t \in A^*$. Then in the case $a < b$ we get $x < y$.

## 2   Ternary Motzkin words

The result of external and internal sorting is an ordered set of words, *ternary Motzkin words*. To get such a set, it is enough to iterate over positive ternary integers starting from 0, checking the balance of 1's and 2's. Again, we consider ternary Motzkin words as character strings in the alphabet $A = \{0, 1, 2\}$; strings are sorted according to alphabetical order $0 < 1 < 2$.

Appendix 1 contains words of nine ranges, 400 codes. Here is the beginning of the ordered set of ternary words:

$M = \{$0, 12, 102, 120, 1002, 1020, 1122, 1200, 1212, 10002, 10020, 10122, 10200, 10212, 11022, 11202, 11220, 12000, 12012, 12102, 12120, 100002, 100020, …$\}$

A Motzkin word cannot end with 1, it can end with either 0 or 2. The first words in the ranges are highlighted in red. In the set $M = \{w_i\}$, $i \in \mathbb{N}$, words are indexed from zero; we will write like this: ind 0 = 0, ind 12 = 1, ind 102 = 2, etc. In the general case, ind $w_i = i$.

In the set of natural numbers ℕ, for each element, there is a single element following it. Let's define the successor function as follows: $\text{Succ}_\mathbb{N}(n) = n+1$, $n \in \mathbb{N}$. The successor function is defined everywhere, the exception is 0, which does not follow any other number.

In M, one can also define the successor function $\text{Succ}_M$, which assigns to $w_i$ the immediately following element $w_{i+1}$, i.e., $\text{Succ}_M(w_i) = w_{i+1}$. The smallest element $w_0$ does not follow any element.

For M and ℕ, the initial elements coincide: $w_0 = n_0 = 0$. In M, ternary codes are grouped by length: $M_1 = \{0\}$, $M_2 = \{12\}$, $M_3 = \{102, 120\}$, and so on. The first two ranges contain one element each, $\#M_1 = \#M_2 = 1$; the next range contains two elements, $\#M_3 = 2$. We work with real Motzkin words, so $\#M_0 = 0$ or $M_0 = \varnothing$.

Let's count the number of elements in the $n$-range. To do this, we subtract the number of elements with the prefix 0 (i.e., the number of words of length $n-1$) from the total number of all $n$-words. So $\#M_n = M_n - M_{n-1} = U_n$, $n \geq 2$. Let's call $U_n$ (the cardinality of $M_n$) the *difference Motzkin number*. The numbers $U_n$ for $n = 1, 2, \ldots$ form the following sequence:

$$1, 1, 2, 5, 12, 30, 76, 196, 512, 1353, 3610, 9713, 26324, 71799, \ldots$$

In the $n$-range of M, there is a minimal word $\min M_n$ and a maximal word $\max M_n$; let's find formulas for these words.

In the $n$-range, the first $n$-word begins with character 1, ends with character 2, and there are $n-2$ zeros inside $n$-word. The minimal $n$-word is written as follows: $1 \cdot 0^{n-2} \cdot 2$, $n \geq 2$. For example, in the case $n = 2$, we get $1 \cdot 0^{2-2} \cdot 2 = 1 \cdot \varepsilon \cdot 2 = 12$. The index of the minimal $n$-word is $M_{n-1}$. Thus,

(1) $\quad \min M_n = 1 \cdot 0^{n-2} \cdot 2, \ \text{ind} \min M_n = M_{n-1}, \ n \geq 2$

In the maximal $n$-word, two characters 12 are repeated; the number of repetitions is $\lfloor n/2 \rfloor$ ($\lfloor x \rfloor$ – the integer part of the real $x$). If $n$ is odd, the maximum is zero-terminated. So the maximal $n$-word is $(12)^k \cdot 0^{n-2k}$, $k = \lfloor n/2 \rfloor$, $n \geq 2$. The index of the maximal $n$-word is $M_n - 1$. As a result, we get the equality

(2) $\quad \max M_n = (12)^k \, 0^{n-2k}, \ \text{ind} \max M_n = M_n - 1, \ k = \lfloor n/2 \rfloor, \ n \geq 2$.

## 3  Motzkin numbers of the second kind

Later we will need ternary Motzkin words for combinatorial analysis and word search in M. But the corresponding numeric codes, *Motzkin numbers of the second kind*, are interesting in their own way, and here we are working not with character strings, but with natural numbers, both in the decimal numeral system and in the ternary numeral system.

## 3.1. Motzkin series

The set of Motzkin numbers of the second kind, hereinafter a *Motzkin series*, we denote by $\mathbb{M} \subset \mathbb{N}$. A Motzkin number of the second kind written in the ternary numeral system, for definiteness, we will mark at the bottom with a three. For example, $12 \in \mathcal{M}$, but $12_3 \in \mathbb{M}$. In Appendix 2, there are over 600 Motzkin numbers of the 2nd kind. Here is the start of the series:

$$\mathbb{M} = \{0, 5, 11, 15, 29, 33, 44/45, 50, 83, 87, 98/99, 104, 116, 128, 132, 135, 140,$$
$$146, 150, 245, 249, 260/261, 266, 278, 290, 294, 297, 302, 308, 312, 332, \ldots\}$$

The series $\mathbb{M} = \{m_i\}$ is completely ordered, the numbers $m_i$ are indexed from 0, like natural numbers. In $\mathbb{M}$ and $\mathbb{N}$, the first numbers coincide, $m_0 = n_0 = 0$. As before, in the ranges, the initial elements are highlighted in red; with another color we marked neighboring numbers that differ by one. Let's call such numbers *Motzkin twins* (recall that there are twins among primes – neighboring odd numbers). Obviously, the indices of the corresponding elements in $\mathbb{M}$ and $\mathcal{M}$ coincide.

The transformations used resemble Gödel's numbering for mathematical objects in a formal language, since in our case, each Motzkin word without a zero suffix (respectively, each Motzkin path without initial horizontal links) is assigned a certain natural number one-to-one. Thus, we use a Motzkin series to *arithmetize* Motzkin words, Motzkin paths and the corresponding graphs. Let's consider an example.

**Example 3.1.** The figure shows the Motzkin path of length 10, which we borrowed from the book by S. Lando [Lan03]. Recall that the left parenthesis, zero, and the right parenthesis correspond to the ascending, horizontal, and descending links of the Motzkin path. The corresponding parenthesis set is (0 ()) () (0), or the word $1012212102 \in \mathcal{M}_{10}$. It is not difficult to arithmetize this path, i.e. to obtain a Motzkin number of the second kind:

$$1012212102_3 = 3^9 + 0 \times 3^8 + 3^7 + 2 \times 3^6 + 2 \times 3^5 + 3^4 + 2 \times 3^3 + 3^2 + 0 \times 3^1 + 2 \times 3^0$$
$$= 19683 + 2187 + 1458 + 486 + 81 + 54 + 9 + 2 = 23960.$$

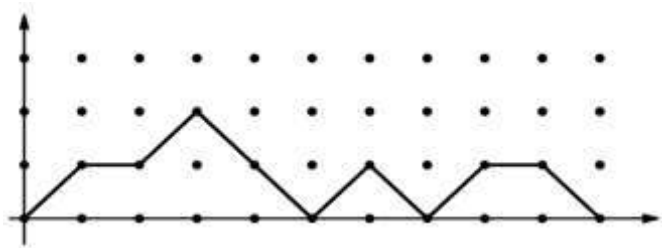

The result is the number $m_{1218} = 23960 \in \mathbb{M}_{10}$ with the index 1218. To perform the reverse procedure, it is enough to recode the decimal number 23960 into the ternary code $1012212102_3$, and then get the parenthesis set and the Motzkin path. □

Let us note some properties of the elements of the set $\mathbb{M}$. For each $m \in \mathbb{M}$ either $m \equiv 0 \pmod{3}$ or $m \equiv 2 \pmod{3}$, since the ternary code cannot end with 1. Multiplying $m \in \mathbb{M}$ by 3 to an arbitrary degree, we get an element of the series again, since each triple shifts the ternary code to the left, adding the suffix 0, which does not change the balance. The tripled number moves into the next range, and the corresponding Motzkin path gets an additional final horizontal link. Thus, the following statement is true.

**Proposition 3.2.** *Let $m \in \mathbb{M}_n$, $n \geq 2$. Then $3^k \times m \in \mathbb{M}_{n+k}$, $k \in \mathbb{N}$.*

In the set $\mathbb{M}$, the successor function moves from a given element $m_i$ to the next element $m_{i+1}$, i.e., $\mathrm{Succ}_{\mathbb{M}}(m_i) = m_{i+1}$. The smallest element $m_0 = 0$ does not follow any element of $\mathbb{M}$.

Let $m_i, m_{i+1} \in \mathbb{M}$ be a pair of twins, then $\mathrm{Succ}_{\mathbb{M}}(m_i) = m_{i+1} = m_i + 1$ and $\mathrm{Succ}_{\mathbb{M}}(m_{i+1}) = m_{i+2} = m_{i+1} + 5$. The second element of the pair is a multiple of 9, so the successor function changes the final two zeros in the ternary code to the suffix $12_3 = 5$. Obviously, $\mathrm{Succ}_{\mathbb{M}}(9m) = 9m + 5$, $m \in \mathbb{M}$.

In the range $\mathbb{M}_n$, there is a minimal element $\min \mathbb{M}_n$ and a maximal element $\max \mathbb{M}_n$. The maximal number of the $n$-range is followed by the minimum of the $(n+1)$-range, that is, $\mathrm{Succ}_{\mathbb{M}}(\max \mathbb{M}_n) = \min \mathbb{M}_{n+1}$. In a range, the maximal number is usually not equal to the minimum, but ranges 1 and 2 each have one element each, so $m_0 = \min \mathbb{M}_1 = \max \mathbb{M}_1 = 0$ and $m_1 = \min \mathbb{M}_2 = \max \mathbb{M}_2 = 5$. Let us consider bounding elements in ranges.

## 3.2. Minimal number in the *n*-range

Let's use the formula (1) to obtain the value of the minimal element in the $n$-range of the set $\mathbb{M}$. In any numeral system, the digits of an integer are numbered from right to left, starting from zero. It is enough for us to sum the least significant 0th digit 2, $2 \times 3^0$, and a 1 in the most significant $(n-1)$-th digit, $1 \times 3^{n-1}$. Thus, the $n$-range starts with the number

$$\min \mathbb{M}_n = 3^{n-1} + 2, \text{ ind } \min \mathbb{M}_n = M_{n-1}, \ n \geq 2,$$

The numbers $\min \mathbb{M}_n$ for $n = 2, 3, \ldots$ form the sequence:

(3)      5, 11, 29, 83, 245, 731, 2189, 6563, 19685, 59051, 177149, …

The infinite sequence (3) is a subsequence in a Motzkin series. The first element is $\min \mathbb{M}_2 = m_1 = 5$, then two equivalent recursive formulas can be used:

$$\min \mathbb{M}_{n+1} = 3 \times (\min \mathbb{M}_n - 2) + 2 = 3 \times (\min \mathbb{M}_n - 1) - 1, \ n \geq 2,$$

or    $\min \mathbb{M}_{n+1} = 4 \times \min \mathbb{M}_n - 3 \times \min \mathbb{M}_{n-1}, \ n \geq 3.$

In the ranges, there are many primes among the minimal numbers.

The sequence (3) is little known, unlike similar Fermat's numbers and Mersenne numbers, whose history is part of the history of Prime numbers (see [Dez11]). However, Neil Sloan's popular On-line Encyclopedia has our series, it's the sequence A168607. The author met the sequence (3) in paper [Mah89], where the representation of integer squares in the ternary number system is considered. In our case, we are interested in special general numbers $E_n = 3^n + 2$, $n \in \mathbb{N}$. Obviously, $\min \mathbb{M}_n = E_{n-1}$.

Note an interesting coincidence: the first five numbers $E_n = 3, 5, 11, 29$, and $83$ are primes, just like the first five Fermat's numbers. Finally, we will talk a little about the numbers $E_n$ and contiguous ('mirror') numbers of the form $\bar{E}_n = 2 \times 3^n + 1$.

## 3.3. Maximal number in the *n*-range

Based on equality (2), let us describe the procedure for calculating the maximal element in the range, starting with the even 4th range. Since $k = \lfloor 4/2 \rfloor = 2$, the transition from a ternary number to a decimal number gives us the value $\max \mathbb{M}_4 = 1212_3 = 12_3 \times 9^0 + 12_3 \times 9^1 = 5 \times (1 + 9^1) = 50$. To get the maximal element in the next odd 5-range, it is enough to triple the obtained sum: $\max \mathbb{M}_5 = 50 \times 3 = 150 = 5 \times 3 \times (1 + 9^1)$ (we simply add 0 to the end of the ternary code). But in order to move from the 5-range to the even 6-range, you need not only to triple the last number, but also to add 5, because in the ternary code, the suffix 00 was formed, which is changed by $12_3 = 5$. The result is $\max \mathbb{M}_6 = 150 \times 3 + 5 = 455 = 5 \times (1 + 9^1 + 9^2)$.

We are dealing with the sum of the first *k* terms of a geometric progression:

(4) $\qquad 1 + 9^1 + 9^2 + \ldots + 9^{k-1} = (9^k - 1)/(9 - 1) = (9^k - 1)/8, \ k \geq 1$.

A formula for the maximal element of the *n*-range takes the form:

(5) $\qquad \max \mathbb{M}_n = 5 \times 3^{n-2k} \times (9^k - 1)/8, \ \text{ind} \max \mathbb{M}_n = M_n - 1, \ k = \lfloor n/2 \rfloor, \ n \geq 2$.

The numbers $\max \mathbb{M}_n$ for $n = 2, 3, \ldots$ form the sequence:

(6) $\qquad$ 5, 15, 50, 150, 455, 1365, 4100, 12300, 36905, 110715, …

The sequence (6) is also an infinite subsequence of a Motzkin series. In (3) and (6) the first elements coincide, because the 2nd range has only one number, which is both a minimum and a maximum. Two recurrence formulas are valid:

$\qquad \max \mathbb{M}_{n+1} = 3 \times \max \mathbb{M}_n$ для $n = 2, 4, 6, \ldots$

and $\quad \max \mathbb{M}_{n+1} = 3 \times \max \mathbb{M}_n + 5$ для $n = 3, 5, 7, \ldots$

All maximal numbers are divided by 5 (the factor of 5 can be a multiple). In sequence (6), the first number is prime, the rest are composite. In odd ranges, the maximal elements are multiples of 3 and 5. The maximum numbers are additionally multiples of 2 if $k = \lfloor n/2 \rfloor$ is even (in this case, the sum (4) has an even number of summands). The resulting formula also works for $n = 1$; since $\mathbb{M}_1 = 0$, the value $\mathbb{M}_2 = 5$ becomes calculated for the second recurrence formula (odd $n$).

## 4  Control points of a Motzkin series

The numbers of the series $\mathbb{M}$ are ordered by length of the ternary codes, the series is divided into ranges: $\mathbb{M} = \mathbb{M}1 \cup \mathbb{M}2 \cup ...$ The range $\mathbb{M}_n$ contains numbers whose ternary codes contain $n$ digits. The first two ranges have one element each; in $n$-range, $n > 2$, we select the minimal element, $n$-minimum, $\min \mathbb{M}_n = E_{n-1} = 3^{n-1} + 2$, and the maximal element, $n$-maximum, $\max \mathbb{M}_n$, which is determined by (5). The element $n$-minimum starts the $n$-range, that has $U_n = M_n - M_{n-1}$ elements, the element $n$-maximum ends the $n$-range. The $n$-minimum index is $M_{n-1}$, the $n$-maximum index is $M_n - 1$. Elements of a Motzkin series are indexed from 0.

It is logical to consider minima (3) and maxima (6) as *control points* of the series (range landmarks), because we know everything about them. The rest of the elements can be defined using the successor function. For example, $\text{Succ}_\mathbb{M}(0) = 5 = 12_3$, $\text{Succ}_\mathbb{M}(5) = 11 = 102_3$, $\text{Succ}_\mathbb{M}(98) = 99$, $\text{Succ}_\mathbb{M}(150) = 245$. In the latter cases, the successor function works with twins $98 = 10122_3$, $99 = 10200_3$ and with elements at the junction of ranges $\max \mathbb{M}_5 = 150 = 12120_3$, $\min \mathbb{M}_6 = E_5 = 245 = 100002_3$.

The successor function for the elements of the Motzkin series is almost identical to the successor function for natural numbers, which maps a natural number $n$ to the number $n+1$ following it. The only difference is the additional check of the balance of 1's and 2's in the ternary code of the subsequent element. For neighboring natural numbers, approximately only in every tenth case, the balance is observed; we called such numbers the twins.

At the junction of the ranges, adjacent elements differ significantly. For example, the jump between ranges 11 and 12 is $\min \mathbb{M}_{12} - \max \mathbb{M}_{11} = 177149 - 110715 = 66434$. Then with each range this value almost triples. In the process of direct enumeration, we can skip ternary numbers starting with 2, but in general case, the direct enumeration is not the best solution to the search procedures.

Additional checkpoints often help us, and there are probably a lot of them. Here are some checkpoints, a few $n$-words (see [Ere19], page 10):

$a^{(n)} = 12 \cdot 0^{n-2}$, ind $a^{(n)} = M_n - M_{n-2}$, $n \geq 2$;
$b^{(n)} = 102 \cdot 0^{n-3}$, ind $b^{(n)} = 2M_{n-1} - M_{n-2} - M_{n-3}$, $n \geq 3$;
$d^{(n)} = 10122 \cdot 0^{n-5}$, ind $d^{(n)} = 2M_{n-1} - M_{n-2} - M_{n-3} - M_{n-5}$, $n \geq 5$.

Usually the index of checkpoints is defined by index equations. For the first checkpoint $a^{(n)}$, the index equality is obvious:

$$\text{ind } a^{(n)} = \text{ind max } \mathbb{M}_n - \text{ind max } \mathbb{M}_{n-2} = (M_n - 1) - (M_{n-2} - 1) = M_n - M_{n-2}.$$

For other points, index equalities were obtained by moving or rearranging the parentheses in Motzkin words. Let's continue Example 3.1.

**Example 4.1.** In Example 3.1, we arithmetized the given Motzkin path and obtained the number $m_{1218} = 23960$ from the 10-range. But how was determined the index 1218? There is a closest checkpoint min $\mathbb{M}_{10} = m_{835} = 3^9 + 2 = 19685$ with index $M_9 = 835$. The next checkpoint max $\mathbb{M}_{10} = 36905$ has the index $M_{10} - 1 = 2187$. We will get the number 1218 faster if we iterate over the numbers from 19685 to 23960 (recall that at each step we need to check the balance of 1's and 2's in the ternary codes of the intermediate numbers). You need to understand that the number of steps during the iteration is not equal to $1218 - 835 = 383$ (which is also a lot), but much more: $23960 - 19685 = 4275$.

We can speed up the search if we connect an additional checkpoint. The checkpoint $b^{(10)} = 102 \cdot 0^{10-3}$ suits us; in this case we get the number $1020000000_3 = 3^9 + 2 \times 3^7 = 19683 + 4374 = 24057$ with index ind $24057 = 2M_9 - M_8 - M_7 = 2 \times 835 - 323 - 127 = 1220$. Now the direct search from the point $m_{1220} = 24057$ down the series is much faster, because the number of steps is $24057 - 23960 = 57$. In fact, we need to "skip" only one number $m_{1219} = 23964$. □

As you can see, the checkpoints in a Motzkin series can speed up enumeration problems, which include finding the next number (when we work with the successor function). However, in large ranges, the distance between the control points is enormous. Let's look at another example.

**Example 4.2.** In the 13th range of a Motzkin series, there are two adjacent elements $m_{25223} = 1021212121212_3 = 686444$ and $m_{25224} = 1100000000022_3 = 708596$ (let the reader believe us). Using the successor function we can write $\text{Succ}_\mathbb{M}(686444) = 708596$. As you can see, to get the next element, you need to jump by the amount of $708596 - 686444 = 22152$. Checkpoints are no help here, and we have to go through more than 22,000 13-digit ternary numbers, checking the balance of 1's and 2's in each.

But there is another simpler and almost instant solution to this problem, namely: we can use the suffix of the ternary code of the original natural number $m_{25223} = 686444$. Working with suffixes, we deal with ternary strings of the set $\mathcal{M}$. In particular, the word $w_{25223} = 1 \cdot 0 \cdot 2 \cdot (12)^5 \in \mathcal{M}_{13}$ has a typical suffix of the form $s = 0 \cdot 2^k \cdot (12)^l$ of length $1 + k + 2l$. In our case $k = 1$ and $l = 5$. The successor function for suffixes, which we will discuss below, gives the suffix $s' = 1 \cdot 0^{2l-1} \cdot 2^{k+1}$ of the same length $1 + k + 2l$ for the next word. The reader will easily get the suffix $s'$ in our case. □

In the next section we consider the successor function, which applies not to the whole word, but only to some part of it, in particular, to the suffix of the ternary Motzkin word. Usually, in any word we can isolate the smallest possible suffix, which, if changed, will give us the next word in the ordered set. We are especially interested in generalized suffixes that work on a subset of an ordered set.

## 5  Successor function for suffixes

In this section we work with character strings, with elements of the set $M$ from the $n$-range, $n \geq 2$, i.e. with Motzkin words that begin with the prefix $p = 1$ and that are bounded by the control points $\min M_n$ and $\max M_n$. Moreover, let us assume that a given $n$-word $w$ is less than the $n$-maximum, and then the next word also belongs to the $n$-range, i.e., $\mathrm{Succ}_M(w) \in M_n$. Let's assume that Motzkin words are checked: the number of 1's and 2's is the same, and in each suffix the balance is negative or zero (there are no more 1's than 2's). Recall that the $k$-digit suffix of a ternary $n$-word is a string of characters of length $k \leq n$.

In a Motzkin word, a digit 1 increases the balance by one, a digit 2 decreases the balance by one, so the balance of the $k$-digit suffix can vary from 0 (zero suffix) to $-k$ (there are only 2's in the suffix). Any suffix can contain leading zeros (unlike a Motzkin word). For example, in the minimal $n$-word $1 \cdot 0^{n-2} \cdot 2$, we are interested in the two-digit suffix 02 (we are not interested in the rest of this word, the prefix $1 \cdot 0^{n-3}$). Hereafter, for brevity, we will sometimes omit the concatenation mark $\cdot$.

**Proposition 5.1.** *Let an $n$-word $w = p \cdot s < \max M_n$ and let the $k$-digit suffix $s \neq w$, i.e., $k < n$ (then the prefix $p \neq \varepsilon$). Then there exists a subsequent $k$-digit suffix $s' > s$ such that* $\mathrm{Succ}_M(w) = p \cdot s'$.

Obviously, for suffix $s$ the *subsequent suffix $s' > s$* has the same length and the same balance. An illustrative example is the suffix 02 with a balance of $-1$. There are two suffixes of the same length 20 and 22 that exceed 02. In terms of balance, the suffix 20 suits us. Therefore, for $s = 02$ we have $s' = 20$. Then $\mathrm{Succ}_M(p \cdot 02) = p \cdot 20$.

A Motzkin word can have several suffixes, and for some of them subsequent suffixes are possible. Obviously, the shortest suffix is to be processed, short suffixes simplify calculations. Thus, the subsequent suffix must be of the same length and balance as the original suffix. A suffix of length $k$ generates *descendants*, $(k+1)$-digit suffixes, if one character is added to the left of the suffix. Below we look at various suffixes, starting with the shortest.

**One-digit suffixes.** Only two one-digit suffixes are possible in a word Motzkin – zero or 2 (the suffix 1 is not allowed, otherwise the balance is upset). The suffixes 0 and 2 are not interchangeable because they have different balances. We consider one-digit suffixes to be at the 1st level of the *suffix tree*. At the 1st level it is impossible to get subsequent suffixes; let's go to the 2nd level of the tree, to the two-digit suffixes.

**Two-digit suffixes.** Five two-digit suffixes are possible in a Motzkin word. Before the final zero can be 0 or 2 (sign 1 is invalid due to imbalance); hence, suffix 0 has two descendants **00** and **20**. The suffix 2 has three descendants – **02**, **12** and **22** (any sign of the alphabet can be in front of the final sign 2). First, let's work with the suffix 00.

**Suffix 00** has a balance of 0, it can only be replaced by suffix 12 with the same balance. Hence, we obtain the first equality for the succession function for suffix: $(00)' = 12$. The remaining suffixes are **20**, **02**, **12**, and **22**.

**Suffix 22** has a balance of $-2$. The sign 2 can be repeated as $2^k$, $k \geq 2$, but there inevitably will be a 0 or 1 in the word before a 2. Thus, the suffix 22 generates two group suffixes: $0 \cdot 2^k$ and $1 \cdot 2^k$, $k \geq 2$. In the suffix $0 \cdot 2^k$, it is sufficient to swap a 0 and the adjacent 2 to get the subsequent suffix $20 \cdot 2^{k-1}$, $k \geq 2$. Recall that we did a similar thing with the suffix 02, and then a similar subsequent suffix 20 was chosen. This means that we can group the suffix $0 \cdot 2^k$ and the suffix 02, if we decrease the lower bound of $k$. The result is equality for the first group suffix:

(7a) $$(0 \cdot 2^k)' = 20 \cdot 2^{k-1}, \ k \geq 1.$$

To get the subsequent suffix for the second group suffix $1 \cdot 2^k$ (length = k+1, balance = $1-k$, $k \geq 2$) just replace 1 with 2, then put two zeros, and we write $k-2$ twos at the end. We get the second subsequent group suffix, which is generated by the suffix 22:

(7b) $$(1 \cdot 2^k)' = 200 \cdot 2^{k-2}, \ k \geq 2.$$

The remaining suffixes are **20** and **12**. Both suffixes go to the next third level of the suffix tree, where we are dealing with three-digit descendants.

**Three-digit suffixes.** At the third level, we have five three-digit suffixes. The suffix 20 has three descendants **020**, **120** and **220**. There are two more descendants **012** and **212** of the suffix 12 (the suffix 112 is invalid due to imbalance).

**Suffix 220** resembles the two-digit suffix 22 (also balance = $-2$). In the case of repeated twos, we get a suffix of the form $2^k 0$ (length = $k+1$, balance = $-k$, $k \geq 2$), to which a 0 or 1 can be added at the front. Thus, suffix 220 generates two group suffixes $0 \cdot 2^k 0$ and $1 \cdot 2^k 0$, $k \geq 2$. In the suffix $0 \cdot 2^k 0$ (length = $k+2$, balance = $-k$, $k \geq 2$) we replace the first 0 with 1 and the final 0 with 2. The result is the subsequent group suffix $1 \cdot 2^k 2 = 1 \cdot 2^{k+1}$, $k \geq 2$. Note that in the case $k = 1$, we get the subsequent suffix 122 for the three-digit suffix 020. And yet, if we take $k = 0$, then we get the subsequent suffix of 12 for the suffix 00. As a result, the resulting suffix $1 \cdot 2^{k+1}$ is subsequent also for the three-digit suffix 020 (if $k = 1$) and for the two-digit suffix 00 (if $k = 0$). Thus,

(8a) $$(0 \cdot 2^k 0)' = 1 \cdot 2^{k+1},\ k \geq 0.$$

For the 2nd group suffix $1 \cdot 2^k 0$ (length $= k + 2$, balance $= 1 - k$, $k \geq 2$), the following suffix will be constructed as follows: we replace the initial 1 with a 2, put $k - 2$ twos at the end (to maintain balance), fill the remaining three positions with zeros. The result is

(8b) $$(1 \cdot 2^k 0)' = 2000 \cdot 2^{k-2},\ k \geq 2.$$

The remaining suffixes are **120**, **012** and **212**.

    **Suffix 012** has zero balance. The subsequent suffix is 102 (also balance $= 0$). Pay attention, in this case we actually moved from the maximum 12 of the 2-range to the minimum 102 of the 3-range. The other maximal elements have a repeating group of 12 in the suffix (and this is the second kind of repeating characters, previously we only worked with a multiple of 2), so it is logical to consider the group suffix $0 \cdot (12)^l$ with length $2l + 1$. Thus,

(9) $$(0 \cdot (12)^l)' = 1 \cdot 0^{2l-1} \cdot 2,\ l \geq 1.$$

The remaining suffixes are **120** and **212**.

    **Suffix 212** (balance $= -1$) has three descendants **0212**, **1212** and **2212**. The first child 0212 can be processed by combining formulas (7a) and (9). For suffix $0 \cdot 2^k (12)^l$ (length $= 1 + k + 2l$, balance $= -k$), we obtain the equality

(10a) $$(0 \cdot 2^k (12)^l)' = 1 \cdot 0^{2l-1} \cdot 2^{k+1},\ k \geq 0,\ l \geq 1.$$

The dependence (10a) absorbs (9). The descendants of suffixes 1212 and 2212 can be treated by formula (10a) if the code starts with 0. Otherwise (if the code starts with 1) the following formula works:

(10b) $$(1 \cdot 2^k (12)^l)' = 2 \cdot 0^{2l+2} \cdot 2^{k-2},\ k \geq 2,\ l \geq 1.$$

    **Suffix 120** (balance $= 0$) has two children 0120 and 2120, which differ from the previous suffixes by a finite zero. The latter circumstance is taken into account when tracing the balance in the subsequent suffixes. The descendants of the suffixes 0120 and 2120 are treated by the following equalities:

(11a) $$(0 \cdot 2^k (12)^l 0)' = 1 \cdot 0^{2l} \cdot 2^{k+1},\ k \geq 0,\ l \geq 1;$$
(11b) $$(1 \cdot 2^k (12)^l 0)' = 2 \cdot 0^{2l+3} \cdot 2^{k-2},\ k \geq 2,\ l \geq 1.$$

This concludes our analysis of the suffixes of the ternary Motzkin words. As the suffix length increased, we got more complex dependencies, some of which absorbed previous, simpler ones. The equations (7a), (10a), (10b), (11a), and (11b) allow us to calculate the subsequent word for an arbitrary element of the ordered set M (in some cases, it may be necessary to correct the lower bounds of the variables $k$ and $l$). We have placed the final equalities in Table 1 (the first column shows the start and last characters of the initial suffixes).

Table 1.

| Type | Initial suffix | Subsequent suffix | Length | Balance | Comments |
|---|---|---|---|---|---|
| 0··2 | $0 \cdot 2^k$ | $20 \cdot 2^{k-1}$ | $k+1$ | $-k$ | $k>0, l=0$ |
| 0··2 | $0 \cdot 2^k (12)^l$ | $1 \cdot 0^{2l-1} \cdot 2^{k+1}$ | $k+2l+1$ | $-k$ | $l>0$ |
| 1··2 | $1 \cdot 2^k (12)^l$ | $2 \cdot 0^{2l+2} \cdot 2^{k-2}$ | $k+2l+1$ | $1-k$ | $k>1$ |
| 0··0 | $0 \cdot 2^k (12)^l 0$ | $1 \cdot 0^{2l} \cdot 2^{k+1}$ | $k+2l+2$ | $-k$ | $k+l>0$ |
| 1··0 | $1 \cdot 2^k (12)^l 0$ | $2 \cdot 0^{2l+3} \cdot 2^{k-2}$ | $k+2l+2$ | $1-k$ | $k>1$ |

The reader can check formulas in Table 1 on the data array in Appendix 1. Let's look at another example.

**Example 5.2.** We have selected a group of the Motzkin words for which subsequent words are obtained. At the beginning of each line, we indicated the index (ordinal number) of the element in M. For each word, the suffix type, length, and balance are shown; we also clarified the values of the variables $k$, $l$.

58. 100·(1220)′ = 1002000, type of suffix 1··0, length = 4, balance = –1, $k = 2$, $l = 0$.
156. 10·(021212)′ = 10100022, type of suffix 0··2, length = 6, balance = –1, $k = 1$, $l = 2$.
241. 11·(122212)′ = 11200002, type of suffix 1··2, length = 6, balance = –2, $k = 3$, $l = 1$.
728. 120·(012120)′ = 120100002, type of suffix 0··0, length = 6, balance = 0, $k = 0$, $l = 2$.
1430. 11011·(02222)′ = 1101120222, type of suffix 0··2, length = 5, balance = –4, $k = 4$, $l = 0$.
1668. 11·(12221212)′ = 1120000002, type of suffix 1··2, length = 8, balance = –2, $k = 3$, $l = 2$.
2496. 1001·(1222120)′ = 10012000002, type of suffix 1··0, length = 7, balance = –2, $k = 3$, $l = 1$.
□

In conclusion, let us make a brief excursion into the special numbers that the author met when analyzing Motzkin numbers of the 2nd kind.

# 6  Special numbers of the form $3^n+2$ and $2 \times 3^n+1$

When encoding parenthesis strings, we replaced the left parenthesis with 1 and the right parenthesis with 2. As a result, the formula for the minimal numbers in ranges of the Motzkin series gives special numbers of the form $E_n = 3^n + 2$ (see [A168607](A168607)). But the reverse coding is also possible: the left parenthesis changes to 2 and the right parenthesis changes to 1. And then a formula for the minimum num-bers in ranges takes the *mirror* form $\bar{E}_n = 2 \times 3^n + 1$ (see A052919). The sequences $E_n$ and $\bar{E}_n$ are practically not described in the literature.

Table 2 shows the initial numbers $E_n$ and $\bar{E}_n$ (primes are highlighted in red). These numbers resemble the famous special *Fermat's numbers* and *Mersenne numbers*, whose history goes back many centuries (Mersenne numbers have been known since ancient times, since the time of Euclid). In the last row, we showed Mersenne numbers $M_n = 2^n - 1$, $n = 1, 2, \ldots$, for comparison (see A000225 ).

Table 2.

| $n$ | 0 | 1 | 2 | 3 | 4 | 5 | 6 | 7 | 8 | 9 | 10 |
|---|---|---|---|---|---|---|---|---|---|---|---|
| $E_n$ | 3 | 5 | 11 | 29 | 83 | 245 | 731 | 2189 | 6563 | 19685 | 59051 |
| $\bar{E}_n$ | 3 | 7 | 19 | 55 | 163 | 487 | 1459 | 4375 | 13123 | 39367 | 118099 |
| $M_n$ | – | 1 | 3 | 7 | 15 | 31 | 63 | 127 | 255 | 511 | 1023 |

As you can see, Table 2 is two-thirds filled with primes $E_n$ and $\bar{E}_n$. The Mersenne primes are much smaller, and this ratio can be traced further, despite the fact that numbers $E_n$ and $\bar{E}_n$ grow much faster than Mersenne numbers. Perhaps in time, it will be possible to prove that the sequences $E_n$ and $\bar{E}_n$ have an infinite number of primes (relative to Mersenne primes, this is still unknown). Let's look at some properties of the numbers $E_n = 3^n + 2$.

All numbers $E_n$ are odd, so $E_n \equiv 1 \pmod{2}$ and in addition $E_n \equiv 2 \pmod{3}$, $n \geq 1$. Similar to Fourier's numbers, a number $E_n$ is not equal to the sum of prime numbers $p + q$, except for $E_1 = 2 + 3$. Let's assume the reverse, $E_n = p + q$, $n > 1$. Since $E_n$ is odd, one of the prime numbers is even, let $q = 2$. Then we obtain $3^n = p$, which is wrong for $n > 1$ and $p > 2$.

Each number $E_n \not\equiv 2 \pmod{p}$, $p \neq 3$. Suppose the opposite, $3^n + 2 \equiv 2 \pmod{p}$ for $p \neq 3$. In this case $3^n \equiv 0 \pmod{p}$, and then $p = 3$. We got a contradiction again.

Let us show that, similarly to Fourier numbers, the number $E_n$ is not a triangular number, except for $E_0$. In other words, $E_n \neq k(k+1)/2$, $k \in \mathbb{N}$ and $n > 0$. The case $k = 2$ is the only exception: $2(2+1)/2 = 3 = E_0$. Indeed, for $k \equiv 0, 2 \pmod{3}$, we obtain $k(k+1)/2 \equiv 0 \pmod{3}$. It remains for us to check option $k \equiv 1 \pmod{3}$. Let us denote $a = k - 1 \equiv 0 \pmod{3}$, then

$$(a+1)(a+2)/2 = 1 + a(a+3)/2 \equiv 1 \pmod{3}.$$

Obviously, the numbers $a$, $a + 3$ have different parity. Again we have a contradiction, since $E_n \equiv 2 \pmod{3}$, $n \geq 1$.

Here's another property: $E_n$ cannot be the square or the cube of an arbitrary natural number. A stronger statement is true, which we will not prove now.

**Proposition 6.1.** $E_n = 3^n + 2 \neq k^s$, $n, k, s \in \mathbb{N}$ and $s > 1$.

We assume that the numbers $E_n$ (as well as the numbers $\bar{E}_n$) have many properties. The sum $E_n + \bar{E}_n = 3(3^n + 1)$ and the difference $\bar{E}_n - E_n = 3^n - 1$ are also interesting in their own way. Let's talk some more about prime numbers.

Table 2 shows the first seven prime numbers $E_n$ (highlighted in red). Here are the indices of the following prime numbers (see here): 14, 15, 24, 36, 63, 98, 110, 112, 123, 126, 139, … Let us show the 17th prime number $E_{126} = 3^{126} + 3 =$ 1310020508637620352391208095712502073964245732475093456566331. This prime is 61 digits long. The closest 12th Mersenne prime by index contains 39 digits: $M_{127} = 170141183460469231731687303715884105727$.

Currently, the 51st Mersenne prime $M_{82589933} = 2^{82589933} - 1$ was found, this prime number contains 24862048 digits. It would be interesting to get something similar for the prime numbers $E_n$.

Gzhel State University, Moscow, 140155, Russia
http://www.art-gzhel.ru/

**Appendix 1.** Motzkin series, ternary words.

*0*: 0, 12, 102, 120, 1002, 1020, 1122, 1200, 1212, 10002, 10020, 10122, 10200;
*13*: 10212, 11022, 11202, 11220, 12000, 12012, 12102, 12120, 100002, 100020, 100122;
*24*: 100200, 100212, 101022, 101202, 101220, 102000, 102012, 102102, 102120, 110022;
*34*: 110202, 110220, 111222, 112002, 112020, 112122, 112200, 112212, 120000, 120012;
*44*: 120102, 120120, 121002, 121020, 121122, 121200, 121212, 1000002, 1000020;
*53*: 1000122, 1000200, 1000212, 1001022, 1001202, 1001220, 1002000, 1002012, 1002102;
*62*: 1002120, 1010022, 1010202, 1010220, 1011222, 1012002, 1012020, 1012122, 1012200;
*71*: 1012212, 1020000, 1020012, 1020102, 1020120, 1021002, 1021020, 1021122, 1021200;
*80*: 1021212, 1100022, 1100202, 1100220, 1101222, 1102002, 1102020, 1102122, 1102200;
*89*: 1102212, 1110222, 1112022, 1112202, 1112220, 1120002, 1120020, 1120122, 1120200;
*98*: 1120212, 1121022, 1121202, 1121220, 1122000, 1122012, 1122102, 1122120, 1200000;
*107*: 1200012, 1200102, 1200120, 1201002, 1201020, 1201122, 1201200, 1201212, 1210002;
*116*: 1210020, 1210122, 1210200, 1210212, 1211022, 1211202, 1211220, 1212000, 1212012;
*125*: 1212102, 1212120, 10000002, 10000020, 10000122, 10000200, 10000212, 10001022;
*133*: 10001202, 10001220, 10002000, 10002012, 10002102, 10002120, 10010022, 10010202;
*141*: 10010220, 10011222, 10012002, 10012020, 10012122, 10012200, 10012212, 10020000;
*149*: 10020012, 10020102, 10020120, 10021002, 10021020, 10021122, 10021200, 10021212;
*157*: 10100022, 10100202, 10100220, 10101222, 10102002, 10102020, 10102122, 10102200;
*165*: 10102212, 10110222, 10112022, 10112202, 10112220, 10120002, 10120020, 10120122;
*173*: 10120200, 10120212, 10121022, 10121202, 10121220, 10122000, 10122012, 10122102;
*181*: 10122120, 10200000, 10200012, 10200102, 10200120, 10201002, 10201020, 10201122;
*189*: 10201200, 10201212, 10210002, 10210020, 10210122, 10210200, 10210212, 10211022;
*197*: 10211202, 10211220, 10212000, 10212012, 10212102, 10212120, 11000022, 11000202;
*205*: 11000220, 11001222, 11002002, 11002020, 11002122, 11002200, 11002212, 11010222;
*213*: 11012022, 11012202, 11012220, 11020002, 11020020, 11020122, 11020200, 11020212;
*221*: 11021022, 11021202, 11021220, 11022000, 11022012, 11022102, 11022120, 11100222;
*229*: 11102022, 11102202, 11102220, 11112222, 11120022, 11120202, 11120220, 11121222;
*237*: 11122002, 11122020, 11122122, 11122200, 11122212, 11200002, 11200020, 11200122;
*245*: 11200200, 11200212, 11201022, 11201202, 11201220, 11202000, 11202012, 11202102;
*253*: 11202120, 11210022, 11210202, 11210220, 11211222, 11212002, 11212020, 11212122;
*261*: 11212200, 11212212, 11220000, 11220012, 11220102, 11220120, 11221002, 11221020;
*269*: 11221122, 11221200, 11221212, 12000000, 12000012, 12000102, 12000120, 12001002;
*277*: 12001020, 12001122, 12001200, 12001212, 12010002, 12010020, 12010122, 12010200;
*285*: 12010212, 12011022, 12011202, 12011220, 12012000, 12012012, 12012102, 12012120;
*293*: 12100002, 12100020, 12100122, 12100200, 12100212, 12101022, 12101202, 12101220;
*301*: 12102000, 12102012, 12102102, 12102120, 12110022, 12110202, 12110220, 12111222;
*309*: 12112002, 12112020, 12112122, 12112200, 12112212, 12120000, 12120012, 12120102;
*317*: 12120120, 12121002, 12121020, 12121122, 12121200, 12121212, 100000002, 100000020;
*325*: 100000122, 100000200, 100000212, 100001022, 100001202, 100001220, 100002000;
*332*: 100002012, 100002102, 100002120, 100010022, 100010202, 100010220, 100011222;
*339*: 100012002, 100012020, 100012122, 100012200, 100012212, 100020000, 100020012;
*346*: 100020102, 100020120, 100021002, 100021020, 100021122, 100021200, 100021212;
*353*: 100100022, 100100202, 100100220, 100101222, 100102002, 100102020, 100102122;
*360*: 100102200, 100102212, 100110222, 100112022, 100112202, 100112220, 100120002;
*367*: 100120020, 100120122, 100120200, 100120212, 100121022, 100121202, 100121220;
*374*: 100122000, 100122012, 100122102, 100122120, 100200000, 100200012, 100200102;
*381*: 100200120, 100201002, 100201020, 100201122, 100201200, 100201212, 100210002;
*388*: 100210020, 100210122, 100210200, 100210212, 100211022, 100211202, 100211220;
*395*: 100212000, 100212012, 100212102, 100212120, 101000022, 101000202, 101000202, ...

**Appendix 2.** Motzkin Number of the Second Kind.

*0*: 0, 5, 11, 15, 29, 33, 44/45, 50, 83, 87, 98/99, 104, 116, 128, 132, 135, 140;
*19*: 146, 150, 245, 249, 260/261, 266, 278, 290, 294, 297, 302, 308, 312, 332, 344;
*35*: 348, 377, 380, 384, 395/396, 401, 405, 410, 416, 420, 434, 438, 449/450, 455;
*51*: 731, 735, 746/747, 752, 764, 776, 780, 783, 788, 794, 798, 818, 830, 834, 863;
*67*: 866, 870, 881/882, 887, 891, 896, 902, 906, 920, 924, 935/936, 941, 980, 992;
*83*: 996, 1025, 1028, 1032, 1043/1044, 1049, 1079, 1115, 1127, 1131, 1136, 1140;
*96*: 1151/1152, 1157, 1169, 1181, 1185, 1188, 1193, 1199, 1203, 1215, 1220, 1226;
*109*: 1230, 1244, 1248, 1259/1260, 1265, 1298, 1302, 1313/1314, 1319, 1331, 1343;
*122*: 1347, 1350, 1355, 1361, 1365, 2189, 2193, 2204/2205, 2210, 2222, 2234, 2238;
*135*: 2241, 2246, 2252, 2256, 2276, 2288, 2292, 2321, 2324, 2328, 2339/2340, 2345;
*148*: 2349, 2354, 2360, 2364, 2378, 2382, 2393/2394, 2399, 2438, 2450, 2454, 2483;
*161*: 2486, 2490, 2501/2502, 2507, 2537, 2573, 2585, 2589, 2594, 2598, 2609/2610;
*174*: 2615, 2627, 2639, 2643, 2646, 2651, 2657, 2661, 2673, 2678, 2684, 2688, 270;
*187*: 2706, 2717/2718, 2723, 2756, 2760, 2771/2772, 2777, 2789, 2801, 2805, 2808;
*200*: 2813, 2819, 2823, 2924, 2936, 2940, 2969, 2972, 2976, 2987/2988, 2993, 3023;
*213*: 3059, 3071, 3075, 3080, 3084, 3095/3096, 3101, 3113, 3125, 3129, 3132, 3137;
*226*: 3143, 3147, 3185, 3221, 3233, 3237, 3320, 3329, 3341, 3345, 3374, 3377, 3381;
*239*: 3392/3393, 3398, 3404, 3408, 3419/3420, 3425, 3437, 3449, 3453, 3456, 3461;
*252*: 3467, 3471, 3491, 3503, 3507, 3536, 3539, 3543, 3554/3555, 3560, 3564, 3569;
*265*: 3575, 3579, 3593, 3597, 3608/3609, 3614, 3645, 3650, 3656, 3660, 3674, 3678;
*278*: 3689, 3690, 3695, 3728, 3732, 3743, 3744, 3749, 3761, 3773, 3777, 3780, 3785;
*291*: 3791, 3795, 3890, 3894, 3905, 3906, 3911, 3923, 3935, 3939, 3942, 3947, 3953;
*304*: 3957, 3977, 3989, 3993, 4022, 4025, 4029, 4040/4041, 4046, 4050, 4055, 4061;
*317*: 4065, 4079, 4083, 4094/4095, 4100, 6563, 6567, 6578/6579, 6584, 6596, 6608;
*330*: 6612, 6615, 6620, 6626, 6630, 6650, 6662, 6666, 6695, 6698, 6702, 6713/6714;
*343*: 6719, 6723, 6728, 6734, 6738, 6752, 6756, 6767/6768, 6773, 6812, 6824, 6828;
*356*: 6857, 6860, 6864, 6875/6876, 6881, 6911, 6947, 6959, 6963, 6968, 6972, 6983/
*369*: 6984, 6989, 7001, 7013, 7017, 7020, 7025, 7031, 7035, 7047, 7052, 7058, 7062;
*382*: 7076, 7080, 7091/7092, 7097, 7130, 7134, 7145/7146, 7151, 7163, 7175, 7179;
*395*: 7182, 7187, 7193, 7197, 7298, 7310, 7314, 7343, 7346, 7350, 7361/7362, 7367;
*408*: 7397, 7433, 7445, 7449, 7454, 7458, 7469/7470, 7475, 7487, 7499, 7503, 7506;
*421*: 7511, 7517, 7521, 7559, 7595, 7607, 7611, 7694, 7703, 7715, 7719, 7748, 7751;
*434*: 7755, 7766/7767, 7772, 7778, 7782, 7793/7794, 7799, 7811, 7823, 7827, 7830;
*447*: 7835, 7841, 7845, 7865, 7877, 7881, 7910, 7913, 7917, 7928/7929, 7934, 7938;
*460*: 7943, 7949, 7953, 7967, 7971, 7982/7983, 7988, 8019, 8024, 8030, 8034, 8048;
*473*: 8052, 8063/8064, 8069, 8102, 8106, 8117/8118, 8123, 8135, 8147, 8151, 8154;
*486*: 8159, 8165, 8169, 8264, 8268, 8279/8280, 8285, 8297, 8309, 8313, 8316, 8321;
*499*: 8327, 8331, 8351, 8363, 8367, 8396, 8399, 8403, 8414/8415, 8420, 8424, 8429;
*512*: 8435, 8439, 8453, 8457, 8468/8469, 8474, 8756, 8768, 8772, 8801, 8804, 8808;
*525*: 8819/8820, 8825, 8855, 8891, 8903, 8907, 8912, 8916, 8927/8928, 8933, 8945;
*538*: 8957, 8961, 8964, 8969, 8975, 8979, 9017, 9053, 9065, 9069, 9152, 9161, 9173;
*551*: 9177, 9206, 9209, 9213, 9224/9225, 9230, 9236, 9240, 9251/9252, 9257, 9269;
*564*: 9281, 9285, 9288, 9293, 9299, 9303, 9323, 9335, 9339, 9368, 9371, 9375, 9386/
*577*: 9387, 9392, 9396, 9401, 9407, 9411, 9425, 9429, 9440/9441, 9446, 9503, 9539;
*590*: 9551, 9555, 9638, 9647, 9659, 9663, 9692, 9695, 9699, 9710/9711, 9716, 9800;
*603*: 9908, 9944, 9956, 9960, 9971, 9983, 9987, 10016, 10019, 10023, 10034/10035;
*615*: 10040, 10070, 10106, 10118, 10122, 10127, 10131, 10142/10143, 10148, 10160;
*626*: 10172, 10176, 10179, 10184, 10190, 10194, 10208, 10212, 10223/10224, 10229;
*637*: 10241, 10253, 10257, 10260, 10265, 10271, 10275, 10295, 10307, 10311, 10340;
*648*: 10343, 10347, 10358/10359, 10364, 10368, 10373, 10379, 10383, 10397, 10401, …